\documentclass[11pt,twoside,a4paper]{article}
\usepackage{amsmath,amssymb,amsthm}

\newtheorem{Thm}{Theorem}[section]
\newtheorem{Lem}[Thm]{Lemma}

\newtheorem{Cor}[Thm]{Corollary}

\theoremstyle{definition}

\theoremstyle{remark}
\newtheorem{Rem}[Thm]{Remark}

\newcommand{\R}{\mathbb{R}}

\newcommand{\N}{\mathbb{N}}

\newcommand{\cA}{\mathcal{A}}
\newcommand{\cB}{\mathcal{B}}
\newcommand{\cC}{\mathcal{C}}

\newcommand{\cN}{\mathcal{N}}
\newcommand{\cO}{\mathcal{O}}

\newcommand{\al}{\alpha}

\newcommand{\ga}{\gamma}

\newcommand{\de}{\delta}
\newcommand{\De}{\Delta}
\newcommand{\ep}{\varepsilon}

\newcommand{\si}{\sigma}

\renewcommand{\phi}{\varphi}

\newcommand{\rank}{\operatorname{rank}}
\newcommand{\hyprank}{\operatorname{rank}_h}
\newcommand{\subcorank}{\operatorname{corank}}
\newcommand{\dist}{\operatorname{dist}}
\newcommand{\diam}{\operatorname{diam}}

\newcommand{\CAT}{\operatorname{CAT}}
\newcommand{\hyp}{\operatorname{H}}
\newcommand{\id}{\operatorname{id}}

\newcommand{\pr}{\operatorname{pr}}
\newcommand{\size}{\operatorname{size}}
\newcommand{\Area}{\operatorname{Area}}

\newcommand{\es}{\emptyset}
\renewcommand{\d}{\partial}
\newcommand{\di}{\d_{\infty}}
\newcommand{\din}{\operatorname{d}_{\infty}}
\newcommand{\set}[2]{\{#1:\,\text{#2}\}}
\newcommand{\sm}{\setminus}
\newcommand{\sub}{\subset}
\newcommand{\sups}{\supset}

\newcommand{\ov}{\overline}
\newcommand{\wt}{\widetilde}

\begin{document}

\title{Hyperbolic rank and subexponential corank of metric spaces}
\author{Sergei Buyalo\footnote{Supported by RFFI Grants 99-01-00104, 
00-15-96024 and SNF Grant 20-57151.99}
\ \& Viktor Schroeder}

\date{}

\maketitle

\begin{abstract} We introduce a new quasi-isometry invariant
$\subcorank X$
of a metric space
$X$
called {\it subexponential corank}.
A metric space
$X$
has subexponential corank
$k$
if roughly speaking there exists a continuous map
$g:X\to T$
such that for each
$t\in T$
the set
$g^{-1}(t)$
has subexponential growth rate in
$X$
and the topological dimension
$\dim T=k$
is minimal among all such maps. Our main result is the inequality
$\hyprank X\le\subcorank X$
for a large class of metric spaces
$X$
including all locally compact Hadamard spaces,
where
$\hyprank X$
is maximal topological dimension of
$\di Y$
among all 
$\CAT(-1)$
spaces
$Y$
quasi-isometrically embedded into
$X$
(the notion introduced by M.~Gromov in a slightly stronger form). This
proves several properties of
$\hyprank$
conjectured by M.~Gromov, in particular, that any Riemannian 
symmetric space
$X$
of noncompact type possesses no quasi-isometric embedding
$\hyp^n\to X$
of the standard hyperbolic space
$\hyp^n$
with
$n-1>\dim X-\rank X$.

\end{abstract}

\section{Introduction}

Given a metric space
$X$
consider all locally compact
$\CAT(-1)$
Hadamard spaces
$Y$
quasi-isometrically embedded into
$X$
and let
$$\hyprank X=\sup_Y\dim\di Y$$
over all such
$Y$.
This quasi-isometry invariant is introduced in \cite[6.$\text{\rm B}_2$]{Gr} 
in a slightly stronger form requiring only
$Y$
to be a geodesic hyperbolic space, where it is called {\it 
hyperbolic corank}. We prefer to call it {\it hyperbolic rank}
reserving the term {\it corank} for a dual notion which is 
central for this paper. We have
$\hyprank X=\dim\di X$
for each
$\CAT(-1)$
space
$X$
(this easily follows from the Morse' quasi-isometry lemma) and
it is known that the Cartesian product 
$X=X_1\times\dots\times X_k$
of Hadamard manifolds with pinched negative curvature,
$-k^2\le K\le -1$,
always has
$$\hyprank X\ge\sum_{i=1}^k\hyprank X_i$$
(see \cite{Gr, BF} for
$X_i$
real hyperbolic manifolds and \cite{FS} for the general case).
Next, let
$X$
be a Riemannian symmetric space of noncompact type. Then
$$\hyprank X\ge\dim X-\rank X$$
(see \cite{Le}). It is conjectured in \cite{Gr} that there is
equality in both cases.

We prove these conjectures by introducing a new quasi-isometry
invariant 
$\subcorank X$
of a metric space
$X$
called {\it subexponential corank} of
$X$
and by showing that
$\hyprank X\le\subcorank X$
for each ``reasonable'' 
$X$.
A metric space
$X$
has subexponential corank
$k$, $\subcorank X=k$,
if roughly speaking there exists a continuous map
$g:X\to T$
such that for each
$t\in T$
the set
$g^{-1}(t)$
has subexponential growth rate in
$X$
and the topological dimension
$\dim T=k$
is minimal among all such maps (for the precise definition
see sect.~\ref{Sect:subcorank}). In this case, we say that 
$X$
supports {\it a subexponential partition of rank}
$k$.
One easily sees that
$\subcorank(X_1\times X_2)\le\subcorank X_1+\subcorank X_2$,
that
$\subcorank X\le\dim X-1$
for every Hadamard manifold
$X$
and that
$\subcorank X\le\dim X-\rank X$
for every Riemannian symmetric space of noncompact type. 

The following theorem, which is our main result, immediately implies 
the above conjectures. Moreover, it works in both directions 
establishing obstructions as to quasi-isometric embeddings of
$\CAT(-1)$
spaces into
$X$
as well as to the existence of subexponential partitions of
$X$. 

We say that a metric space
$X$
has the {\it QPC property} or 
$X$
is a QPC-space if every quasi-isometric map
$f:Y\to X$
is parallel to a continuous one, i.e., there exists a continuous
$f':Y\to X$
with
$\dist(f'(y),f(y))\le C<\infty$
for all
$y\in Y$
(QPC stays for Quasi-isometric maps are Parallel to Continuous ones).
For example, every locally compact Hadamard space is QPC (for more 
details see sect.~\ref{Sect:subcorank}).
 
\begin{Thm}\label{Thm:mainresult} Let
$X$
be a metric space which is quasi-isometric to a QPC one. Then
$$\hyprank X\le\subcorank X.$$
In particular, there is no quasi-isometric embedding
$Y\to X$
of 
$\CAT(-1)$
spaces
$Y$
with
$\dim\di Y>\subcorank X$
and if 
$X$
is QPC then it supports no subexponential partition of rank less than
$\hyprank X$.

\end{Thm}

From \cite{BF, FS, Le} and Theorem~\ref{Thm:mainresult}
we obtain 

\begin{Cor}\label{Cor:equality} (1) Let 
$X$ be the Cartesian product of Hadamard manifolds with
pinched negative curvature,
$X=X_1\times\dots\times X_k$.
Then
$$\hyprank X=\sum_{i=1}^k\hyprank X_i=\sum_{i=1}^k\subcorank X_i
  =\subcorank X.$$

(2) Let
$X$
be a Riemannian symmetric space of noncompact type. Then
$$\hyprank X=\dim X-\rank X=\subcorank X.$$

\end{Cor}

\begin{Rem} We use the
$\CAT(-1)$
condition in the definition of
$\hyprank$
to simplify the proof of Theorem~\ref{Thm:mainresult}.
It can be modified in a way to include all (complete, 
locally compact) geodesic hyperbolic spaces proving by that Gromov's 
conjectures in their full generality.
\end{Rem}

The invariant
$\subcorank$
and Theorem~\ref{Thm:mainresult} have further applications.
For instance, we prove

\begin{Thm}\label{Thm:graphmfd} Let
$X$
be the metric universal covering of a closed 3-dimensional
nonpositively curved graph manifold. Then
$\subcorank X\le 1$.
As a consequence,
$Z=X_1\times\dots\times X_n$
has
$\subcorank Z\le n$
and possesses no quasi-isometric embedding
$Y\to Z$
of 
$\CAT(-1)$
spaces
$Y$
with
$\dim\di Y>n$,
where each 
$X_i$
is the universal covering of a 3-dimensional nonpositively
curved graph manifold.

\end{Thm}

Note that each
$X_i$
in Theorem~\ref{Thm:graphmfd} is not a 
$\CAT(-1)$
space. It only
possesses a quasi-isometric
embedding
$\hyp^2\to X_i$,
hence
$Z=X_1\times\dots\times X_n$
has an exponential growth rate and
$\hyprank Z=n$
by combining the results of \cite{BF} and Theorem~\ref{Thm:graphmfd}.

Another consequence of Theorem~\ref{Thm:mainresult} is that
$\hyprank X=n-1=\subcorank X$
for every
$\CAT(-1)$
Hadamard {\it manifold}
$X$, $\dim X=n$.
One might expect that the equality
$\hyprank X=\subcorank X$
is true for each
$\CAT(-1)$
space
$X$. 
However, this is not the case. For example, let
$X$
be a simplicial tree whose vertices have degree
at least 3 and whose edges have length 1. Then
$X$
is a 
$\CAT(-1)$
space with
$\hyprank X=\dim\di X=0$
while
$\subcorank X=1$
because every continuous map
$X\to T$
into 0-dimensional space
$T$
is constant and 
$X$
has exponential growth rate.

We conclude the Introduction by a sketch of our initial proof that 
there is no quasi-isometric embedding
$f:\hyp^3\to\hyp^2\times\R^n$
for any
$n\ge 0$.
This proof was the starting point to introduce the subexponential
corank and finally to prove Theorem~\ref{Thm:mainresult}. 
Note that any simple counting argument does not work
because both the source space
$\hyp^3$
and the target space
$\hyp^2\times\R^n$
have exponential growth rates.

Assume that there is such an
$f$. Assume for simplicity that
$f$
is smooth and biLipschitz, that is
$\frac{1}{a}\dist(x,x')\le\dist(f(x),f(x'))\le a\dist(x,x')$
for some
$a\ge 1$
and all
$x$, $x'\in \hyp^3$.
Fix a horosphere
$T\sub\hyp^2$
and define a map
$g:\hyp^2\times\R^n\to T$
by projecting onto the first factor and then
projecting onto
$T$
along geodesics in
$\hyp^2$
orthogonal to
$T$.
Note that all its fibers
$g^{-1}(t)$, $t\in T$
are isometric to
$\R^{n+1}$
and geodesically embedded in
$\hyp^2\times\R^n$.
Composing 
$g$
with
$f$
and identifying
$T=\R$
we obtain a (smooth) function
$h:\hyp^3\to\R$, $h=g\circ f$. 

Next, we fix
$y_0\in\hyp^3$
and consider the metric sphere
$S_R\sub\hyp^3$
of a sufficiently large radius
$R$
centered at
$y_0$.
Considered with the induced intrinsic metric, 
$S_R$
is isometric to the standard sphere
$S_{\rho}^2\sub\R^3$
of radius
$\rho$
which is exponentially large in
$R$, $\rho\sim e^R$.
We may further assume for simplicity that
$h$
restricted to
$S_R$
is a Morse function,
$h:S_{\rho}^2\to\R$.
Now we come to the crucial point. Every level set
$\ga_t=h^{-1}(t)$, $t\in\R$
is mapped by
$f$
into the
$(n+1)$-flat
$F=g^{-1}(t)\sub\hyp^2\times\R^n$.
Since
$f$
is Lipschitz and 
$F$
is geodesic,
$f(\ga_t)$
sits in a ball in
$F$
of radius
$\sim R$
with respect to the induced {\it intrinsic} metric on
$F$. 
Thus
$f(\ga_t)$
has a polynomial in
$R$
size, where by ``size'' we mean, for example, the minimal
number of points of some separated net in
$f(\hyp^3)$
needed to cover
$f(\ga_t)$
by balls of a fixed radius centered at these points. Since
$f$
is biLipschitz, this implies that every connected component of
$\ga_t\sub S_{\rho}^2$
(but may be not the whole 
$\ga_t$!)
has the diameter polynomial in
$R$,
i.e., essentially smaller than
$\rho$.
This conclusion is a key point which easily leads to 
a contradiction. Having all connected components of all level
sets
$\ga_t$
sufficiently small compare to
$\rho$
one can, for example, continuously contract the sphere
$S_{\rho}^2$
to a point (with some little work near critical levels of
index one). However, these details are inessential for
the generalization and we omit them.

\medskip
{\bf Acknowledgment.} The first author is happy to express his
deep gratitude to the University of Z\"urich for the support,
hospitality and excellent working conditions while writing the 
paper.

\section{Subexponential corank}\label{Sect:subcorank}

In this section we define the subexponential corank of a metric space
$X$
and establish some of its properties which are needed for the proof 
and applications of Theorem~\ref{Thm:mainresult}.

\subsection{Preliminaries}

The key notions of 
$\CAT(-1)$-spaces
and Hadamard spaces are well established by now, and the
reader may consult, for instance, \cite{BH} for them.
 
Given a metric space
$X$
and
$x\in X$,
we denote by
$B_R(x)$
the open ball of radius
$R$
in
$X$
centered at
$x$
and by
$\ov B_R(x)$
the corresponding closed ball. A subset
$A\sub X$
is called {\it net} if
$\dist(A,X)<\infty$.
$A$
is
$\de$-{\it separated},
$\de>0$,
if $\dist(a,a')\ge\de$
for every distinct
$a$, $a'\in A$.
We usually use the notation 
$X_{\de}$
for a separated net in
$X$,
where
$\de$
is the separation constant of
$X_{\de}$.
Note that the balls 
$B_{\de}(a)$
centered at the points
$a\in X_{\de}$
of some {\it maximal} separated net
$X_{\de}\sub X$
cover
$X$.

Assume that a maximal separated net
$X_{\de}\sub X$
and
$\si\ge\de$
are fixed. Then we define the size of 
$A\sub X$
(w.r.t.
$X_{\de}$
and
$\si$)
as the number
$\size_{X_{\de},\si}(A)\in\N\cup\{\infty\}$
of points
$x\in X_{\de}$
with
$B_{\si}(x)\cap A\not=\es$.
By the remark above, the union of all such balls contains
$A$.

A map
$f:X\to Y$
between metric spaces
$X$, $Y$
is said to be {\it quasi-isometric}, if
$$\frac{1}{a}\dist(x,x')-b\le\dist\left(f(x),f(x')\right)
   \le a\dist(x,x')+b$$
for some
$a\ge 1$, $b\ge 0$
and all
$x$, $x'\in X$.
In this case we say that
$f$
is
$(a,b)$-quasi-isometric.
If in addition
$f(X)$
is a net in
$Y$
then
$f$
is called {\it a quasi-isometry} and the spaces
$X$
and
$Y$
are {\it quasi-isometric}.

Recall that the (topological) dimension of a compact set
$K$
is the minimal integer
$n=\dim K$
such that one can inscribe in every open covering of
$K$
a finite closed covering having the multiplicity
$\le n+1$.
If in addition
$K$
is a metric space then
$\dim K$
is the same as the minimal number
$n$
having the property that for every
$\ep>0$
there is a finite closed covering of
$K$
by sets with diameter
$<\ep$
and the multiplicity
$\le n+1$.
For a topological space
$X$
we use the definition
$\dim X=\sup\set{\dim K}{$K\sub X\ \text{is compact}$}$.

\subsection{Definition of the subexponential corank}

{\it A continuous partition} 
$$X=\bigcup_{t\in T}g^{-1}(t)$$
of a metric space
$X$
is given by a continuous map
$g:X\to T$.
We use the notation
$(X,g,T)$
and say that the partition
$(X,g,T)$
has rank
$k=\sup\set{\dim g(K)}{$K\sub X\ \text{is compact}$}$.
Fix
$x_0\in X$.
A continuous partition
$(X,g,T)$
is said to be {\it subexponential} if the following 
holds. For each maximal separated net
$X_{\de}\sub X$
with a sufficiently large separation constant
$\de$,
each sufficiently large
$\si\ge\de$
and every 
$\ep>0$
there exists
$R_0=R_0(X_{\de},\si,\ep)\ge 1$
such that for every  
$R\ge R_0$
and every
$t\in T$
we have
$$\frac{1}{R}\ln\size_{X_{\de},\si}(g^{-1}(t)\cap\ov B_R(x_0))<\ep.$$
This property, obviously, is independent of the choice of
$x_0\in X$.
Now we define {\it the subexponential
corank} of
$X$
as
$$\subcorank X=\sup_{Z\sim X}\inf\rank(Z,g,T),$$
where the supremum is taken over all
$Z$
quasi-isometric to
$X$
and the infimum is over all subexponential partitions of
$Z$.
Clearly,
$\subcorank X$
is a quasi-isometric invariant.

\begin{Rem} It is the controversy between continuity conditions needed
for the proof of Theorem~\ref{Thm:mainresult} and the quasi-isometry 
invariance needed for its applications which makes the definition of 
$\subcorank X$
rather entangled. In particular, taking the supremum over all
$Z$
quasi-isometric to
$X$
in the definition of
$\subcorank X$
is necessary because for any discrete space
$Z$
the trivial partition
$(Z,\id,Z)$
is subexponential and has rank 0, i.e., minimal possible rank.
\end{Rem}

\begin{Rem} The property of a partition
$(X,g,T)$
to be subexponential means roughly speaking that every fiber
$g^{-1}(t)$
has a subexponential growth rate and a bounded distortion in
$X$.
The last condition is essential. For example, all fibers of 
the partition
$(\hyp^n,g,\R)$
given by a Busemann function
$g:\hyp^n\to\R$
are isometric to
$\R^{n-1}$
in the induced Riemannian metric and hence they have a subexponential 
(in fact, polynomial) growth rate.
However, this partition is by no means subexponential because
each horosphere
$g^{-1}(t)$, $t\in\R$
is exponentially distorted in
$\hyp^n$
and for a fixed
$x_0\in\hyp^n$
the balls
$\ov B_R(x_0)$
contain exponentially large pieces of it.
\end{Rem}

We shall frequently use the following actually straightforward

\begin{Lem}\label{Lem:reduction1} If
$f:X\to Z$
is a {\it continuous} quasi-isometric map and
$(Z,g,T)$
is a subexponential partition then
$(X,g\circ f,T)$
is a subexponential partition too.
\end{Lem}

\begin{proof} Assuming that
$f$
is
$(a,b)$-quasi-isometric,
we fix
$x_0\in X$
and put
$z_0=f(x_0)$.
Let
$\de_0$, $\si_0\ge\de_0$
be the separation and radius constants respectively for
$Z$
involved in the subexponential size estimate. We put
$\de=a(\de_0+b)$
and take a maximal separated net
$X_{\de}\sub X$.
Then
$f(X_{\de})$
is
$\de_0$-separated,
hence, it is contained in some maximal separated net
$Z_{\de_0}\sub Z$.

Next, we fix
$\si\ge\max\{\de,\frac{1}{a}(\si_0-b)\}$
and take
$t\in T$.
If 
$B_{\si}(x)$
intersects
$(g\circ f)^{-1}(t)\cap\ov B_R(x_0)$
for some
$x\in X_{\de}$
then
$f(B_{\si}(x))$ 
intersects
$g^{-1}(t)\cap\ov B_{aR+b}(z_0)$
since
$f(\ov B_R(x_0))\sub\ov B_{aR+b}(z_0)$.
Thus
$B_{a\si+b}(f(x))$
intersects
$g^{-1}(t)\cap\ov B_{aR+b}(z_0)$
and consequently we have
$$\size_{X_{\de},\si}\left((g\circ f)^{-1}(t)\cap\ov B_R(x_0)\right)
   \le\size_{Z_{\de_0},a\si+b}\left(g^{-1}(t)\cap\ov B_{aR+b}(x_0)\right)$$
for every
$R>0$, $t\in T$.   
Fix
$\ep>0$.
Then for
$aR+b\ge R_0(Z_{\de_0},a\si+b,\frac{\ep}{a+b})$
we obtain
$$\frac{1}{aR+b}
  \ln\size_{X_{\de},\si}\left((g\circ f)^{-1}(t)\cap\ov B_R(x_0)\right)
  <\frac{\ep}{a+b}.$$
Thus taking
$R_0(X_{\de},\si,\ep)\ge\max\left\{1,\frac{1}{a}
  \left[R_0\left(Z_{\de_0},a\si+b,\ep/(a+b)\right)-b\right]\right\}$
we obtain the required estimate 
$$\frac{1}{R}\ln\size_{X_{\de},\si}
  \left((g\circ f)^{-1}(t)\cap\ov B_R(x_0)\right)<\ep$$
for every
$R\ge R_0(X_{\de},\si,\ep)$
and
$t\in T$.
\end{proof}

\begin{Lem}~\label{Lem:reduction2} Assume that
a metric space
$X$
is quasi-isometric to a QPC space
$Z$.
Then
$$\subcorank X=\inf\rank(Z,g,T),$$
where the infimum is taken over all subexponential
partitions of
$Z$.
\end{Lem}

\begin{proof} If
$X'$
is quasi-isometric to
$X$
then it is quasi-isometric to
$Z$
too.
Any map
$X'\to Z$
parallel to a quasi-isometric one is quasi-isometric, thus
there is a continuous quasi-isometry
$X'\to Z$.
By Lemma~\ref{Lem:reduction1} we have
$\inf\rank(X',g',T')\le\inf\rank(Z,g,T)$. 
Hence, the claim.
\end{proof}

Now we explain why all locally compact Hadamard spaces
are QPC, though the argument (which is basically standard) 
might be applied to a much broader class of metric spaces.

\begin{Lem}\label{Lem:QPC} Every locally compact Hadamard
space
$X$
is a QPC-space.
\end{Lem}

\begin{proof} Assume that we have an
$(a,b)$-quasi-isometric
map
$f:Y\to X$.
We take a maximal 
$\de$-separated
net
$Y_{\de}\sub Y$
with
$a\de+b\ge\de_0>0$
and note that every ball
$B_{2\de}(\al)$
with
$\al\in Y_{\de}$
contains only a finite number of elements of
$Y_{\de}$.
This is because
$f(Y_{\de})$
is
$\de_0$-separated
and
$X$
is finitely compact, thus the ball
$B_{2a\de+b}\left(f(\al)\right)\sups f\left(B_{2\de}(\al)\right)$
intersects
$f(Y_{\de})$
over a finite set. It follows that the nerve
$\cN$
of the covering
$\cA=\set{B_{\de}(\al)}{$\al\in Y_{\de}$}$
of
$Y$
is a locally finite simplicial complex. Choosing
a continuous partition of unity
$\set{p_{\al}:Y\to\R}{$\al\in Y_{\de}$}$
subordinate to
$\cA$,
we obtain a continuous map
$g:Y\to\cN$
by
$$g(y)=\sum_{\al\in Y_{\de}}p_{\al}(y)\al,$$
where we identify
$Y_{\de}$
with the 0-skeleton of
$\cN$.
Note that
$g(y)$
lies in the simplex 
$\De_y$
spanned by
$\set{\al\in Y_{\de}}{$\dist(\al,y)<\de$}$.

Next, we extend
$f|Y_{\de}:\text{ske}_0\cN\to X$
to a continuous
$\ov f:\cN\to X$
using convexity of
$X$
and acting by the induction on the dimension of the skeletons.
Then
$\ov f(\De_y)\sub B_{a\de+b}(f(y))$
and thus
$\ov f\circ g:Y\to X$
is a continuous map parallel to
$f$.
\end{proof}

\subsection{Properties of the subexponential corank}

We list some properties of
$\subcorank$
which easily follow from the definition.

(1) $\subcorank X\le\dim X$
for every QPC-space
$X$.
This immediately follows from Lemma~\ref{Lem:reduction2}.

(2) $\subcorank(\R^n)=0$
for each
$n\ge 0$.
Moreover, if the volume entropy
$$h(X)=\lim\sup_{R\to 0}\frac{1}{R}\ln\size(B_R(x_0))=0,$$
then
$\subcorank X=0$.
Note that the condition
$h(X)=0$
is a quasi-isometry invariant. In this case, the trivial partition
$(X,g,\{\text{pt}\})$
given by a constant map
$g$
is subexponential.

(3) $\subcorank(X_1\times X_2)\le
   \subcorank(X_1)+\subcorank(X_2)$ 
if 
$X_1\times X_2$
is a QPC-space. In this case, both
$X_1$, $X_2$
are QPC and the product partition
$(X_1\times X_2,g_1\times g_2,T_1\times T_2)$
of subexponential partitions
$(X_i,g_i,T_i)$, $i=1,2$
is subexponential. Finally, it is well known (see
\cite{HW}) that
$\dim(A\times B)\le\dim A+\dim B$
if at least one of the spaces
$A$, $B$
is not empty.

(4) Let
$X$
be a Hadamard manifold. Then
$\subcorank X\le\dim X-1$.
Projecting onto a horosphere
$T\sub X$
along geodesics orthogonal to
$T$,
we obtain the subexponential partition
$(X,g,T)$
of rank 
$\dim X-1$.
Its fibers
$g^{-1}(t)$
are geodesics.

(5) Let
$X$
be a Riemannian symmetric space of noncompact type. Then
$\subcorank X\le\dim X-\rank X$.
An Iwasawa decomposition
$G=NAK$
of the connected component of the identity in its isometry
group allows to identify
$X$
with the solvable group
$NA$.
Fix
$x_0\in X$
and consider the orbit
$T=Nx_0\sub X$.
Then the map
$g:X\to T$
given by
$g(x)=nx_0$, 
where
$x=nax_0$
for
$n\in N$, $a\in A$,
defines the subexponential partition
$(X,g,T)$
of rank
$\dim T=\dim N=\dim X-\rank X$.
Its fibers
$g^{-1}(t)$
are geodesic
$\rank X$-flats.

(6) If
$X$
is quasi-isometric to a QPC-space and
$X'$
admits a quasi-isometric embedding into
$X$
then
$\subcorank X'\le\subcorank X$.
This follows Lemma~\ref{Lem:reduction1} and Lemma~\ref{Lem:reduction2}.

\section{Hyperbolic rank and subexponential corank}\label{Sect:proofmainresult}
Here we prove Theorem~\ref{Thm:mainresult} and its corollaries. The proof of
Theorem~\ref{Thm:mainresult} in outline goes as follows. By property~(6)
above it suffices to show that
$\hyprank Y\le\subcorank Y$
for each locally compact
$\CAT(-1)$
Hadamard space
$Y$.
Since
$\hyprank Y=\dim\di Y$
and
$Y$
is certainly QPC it suffices to show that
$\dim\di Y\le\rank(Y,g,T)$
for any subexponential partition
$(Y,g,T)$.

We have a continuous map
$\pr_R:\di Y\to S_R$,
where
$S_R\sub Y$
is the metric sphere of radius
$R$
centered at some fixed point
$y_0\in Y$,
given by the intersection of the ray
$y_0\xi$, $\xi\in\di Y$
with
$S_R$.
Using the
$\CAT(-1)$ condition
we introduce a metric
$\din$
on
$\di Y$
for which
$\din(\xi,\xi')\sim e^{-\dist(y_0,\xi\xi')}$
for
$\xi$, $\xi'\in\di Y$,
where
$\xi\xi'$
is the geodesic in
$Y$
asymptotic to
$\xi$
in one direction and to
$\xi'$
in the other one. Note that
$\di Y$
is compact since
$Y$
is complete and locally compact and that the metric topology of
$\din$
coincides with the topology of
$\di Y$.
Thus the resulting map
$h_R=g\circ\pr_R:\di Y\to T$
is continuous and its image
$K_R=h_R(\di Y)$
is compact.

Now, we can inscribe in any open covering 
$\cO$
of
$K_R$
a finite closed covering 
$\cC$
having the multiplicity
$\le n+1$, $n=\rank(Y,g,T)$.
Since the partition
$(Y,g,T)$
is subexponential, we may produce for each sufficiently large
$R$
an open covering
$\cO_R$
of
$K_R$
such that the size of the preimage in
$S_R$
of every element of
$\cO_R$
is subexponential in
$R$. 
Consequently, we can find a finite closed covering 
$\cC_R$
of
$K_R$
with the multiplicity
$\le n+1$
such that the size of the preimage in
$S_R$
of every its element 
$C$
is subexponential in
$R$.
However, it does not mean that the diameter of
$g^{-1}(C)$
measured in
$S_R$
is subexponential in
$R$
because it may have different connected components far away
from each other. Nevertheless, we can rearrange the finite
closed covering
$\cA_R=g^{-1}(\cC_R)$
of
$S_R$
in such a way that it is still finite, has the multiplicity
$\le n+1$
and the diameter of every its element is subexponential in
$R$.
In this way, we produce a finite closed covering
$\cB_R=\pr_R^{-1}(\cA_R)$
of
$\di Y$
of multiplicity
$\le n+1$
with maximal diameter of its elements measured in the metric
$\din$
going to 0 as
$R\to\infty$.
Hence,
$\dim\di Y\le n$.
Q.E.D.

We fill in details answering the 
following questions

\begin{itemize}

\item[(1)] how do we choose the open covering
$\cO_R$
of
$K_R$;

\item[(2)] how to rearrange
$\cA_R$;

\item[(3)] how do we define the metric
$\din$;

\item[(4)] how to estimate
$\diam(B)$
w.r.t.
$\din$
for each
$B\in\cB_R$.

\end{itemize}

(1) Fix a maximal
separated net
$Y_{\de}\sub Y$
and
$\si\ge\de$,
which are sufficiently large to fit the
$(Y,g,T)$
conditions,
and a base point
$y_0\in Y$.
The notation
$\size_{Y_{\de},\si}(A)$
will be abbreviated to
$\size(A)$
for
$A\sub Y$.
Fix
$\ep$
with
$0<\ep<1$.
For each
$R>0$
and
$t\in T$,
the set
$g^{-1}(t)\cap S_R(y_0)$
is covered by
$N(R,t)=\size\left(g^{-1}(t)\cap S_R(y_0)\right)$
open balls
$B_{\si}(y)$
with
$y\in Y_{\de}$.
Using properties of the subexponential partition
$(Y,g,T)$
we have 
$$\frac{1}{R}\ln N(R,t)\le
  \frac{1}{R}\ln\size\left(g^{-1}(t)\cap\ov B_R(y_0)\right)<\ep$$
for every sufficiently large
$R\ge 1$
and every
$t\in T$.

Let
$V\sub Y$
be the union of those balls. 
We claim that there is a neighborhood
$U_t$
of
$t$
with
$g^{-1}(U_t)\cap S_R(y_0)\sub V$.
Otherwise, we find a sequence
$y_i\in S_R(y_0)\sm V$
with
$g(y_i)\to t$.
Since 
$Y$
is complete and locally compact, we may assume that
$y_i\to y_{\infty}\in S_R(y_0)\sm V$.
On the other hand,
$y_{\infty}\in g^{-1}(t)$
by continuity of
$g$.
This is a contradiction.

Then
$\cO_R=\set{U_t}{$t\in K_R$}$
is the required open covering of
$K_R$:
we have
$$\frac{1}{R}\ln\size\left(g^{-1}(U_t)\cap S_R(y_0)\right)<\ep$$
for each
$t\in K_R$.

(2) For each
$A\in\cA_R$
we consider its covering by balls
$B_{\si}(y)$, $y\in Y_{\de}$
which intersect
$A$.
The union of these balls has only a finite number of
connected components
$\al$
since every such
$\al$
contains at least one ball
$B_{\si}(y)$
and there is only a finite number of such balls in
$B_{R+\si}(y_0)$.
Now, we decompose
$A$
into the finite disjoint union
$A=\cup_{\al}A_{\al}$,
where each
$A_{\al}$
is the intersection of
$A$
with a connected component
$\al$.
Thus the diameter of
$A_{\al}$
measured in
$S_R$
is at most
$2\de\size(A)$,
$\diam_R(A_{\al})\le 2\de e^{\ep R}$,
and this gives the desired rearrangement of
$\cA_R$
(for which we use the same notation).

(3) We define the metric
$\din$
on
$\di Y$
as follows. Given
$\xi$, $\xi'\in\di Y$,
we consider the unit speed geodesic rays
$c_{\xi}$, $c_{\xi'}$
from
$y_0$
to
$\xi$, $\xi'$
respectively and put
$$\din(\xi,\xi')=\lim_{s\to\infty}
  \angle(\ov{c_{\xi}}(s)\ov y_0\ov{c_{\xi'}}(s)),$$
where
$\angle(\ov{c_{\xi}}(s)\ov y_0\ov{c_{\xi'}}(s))$
is the angle at
$\ov y_0$
of the comparison triangle in
$\hyp^2$
for the triangle
$c_{\xi}(s)y_0c_{\xi'}(s)$.
From the hyperbolic geometry we know that
$$\tan\left(\frac{1}{4}\din(\xi,\xi')\right)=
  e^{-\dist(\ov y_0,\ov\xi\ov\xi')}.$$
Thus
$\din(\xi,\xi')\le 4e^{-\dist(y_0,\xi\xi')}$.  

(4) Given
$B\in\cB_R$
let
$\xi$, $\xi'\in B$
be two remotest points
(w.r.t.
$\din$), $v=\pr_R(\xi)$, $v'=\pr_R(\xi')\in S_R$.
Then
$v$, $v'\in A$
for
$A\in\cA_R$
with
$\pr_R^{-1}(A)=B$.
In particular,
$\diam_R(A)\le 2\de e^{\ep R}$.
We take the midpoint
$y$
of the geodesic segment
$vv'\sub Y$
and put
$\rho=\dist(y_0,y)$, $r=\dist(y,S_R)$.
Then
$\rho+r\ge R$,
and an easy comparison argument shows that the length
of any curve in
$S_R$
between
$v$, $v'$
and hence
$\diam_R(A)$
is at least as large as the length of a half-circle
of radius
$r$
in
$\hyp^2$
(we use here that the points
$v$, $v'$
are visible from
$y$
at the angle
$\pi$).
Thus we obtain
$r\le \ep R+c$
for some
$c=c(\de)$.
It follows 
$$\dist(y_0,\xi\xi')\ge\rho\ge R-r\ge(1-\ep)R-c.$$
Since
$\ep<1$,
we have
$$\diam(B)=\din(\xi,\xi')\le 4e^ce^{-(1-\ep)R}\to 0
  \quad\text{as}\quad R\to\infty.$$
This completes the proof of Theorem~\ref{Thm:mainresult}.  

\begin{proof}[Proof of Corollary~\ref{Cor:equality}] (1) By 
Lemma~\ref{Lem:QPC},
$X$
is QPC. Thus
$\subcorank X\le\sum_{i=1}^k\subcorank X_i\le
  \sum_{i=1}^k(\dim X_i-1)$
by Properties~(3), (4), sect.~\ref{Sect:subcorank}.
On the other hand, it follows from \cite{FS} that
$\hyprank X\ge\sum_{i=1}^k\hyprank X_i=\sum_{i=1}^k
  (\dim X_i-1)$.
By Theorem~\ref{Thm:mainresult},
$\hyprank X\le\subcorank X$.
Hence, the claim.  

(2) We have
$\dim X-\rank X\le\hyprank X$
by \cite{Le},
$\subcorank X\le\dim X-\rank X$
by Property~(5), and
$\hyprank X\le\subcorank X$
by Theorem~\ref{Thm:mainresult}. Hence, the claim.
\end{proof}

\section{An exotic subexponential partition} 

Here we prove Theorem~\ref{Thm:graphmfd}. Recall 
(see, for instance, \cite{BS, CK, KL})
that 
$X$
can be represented as the countable union
$X=\cup_vX_v$
of blocks, where each
$X_v$
is a closed convex subset in
$X$
isometric to the metric product
$F_v\times\R$
and 
$F_v$
is the universal covering of a compact nonpositively
curved surface with geodesic boundary. Every two blocks
are either disjoints or intersect over a boundary
component which is a 2-flat in
$X$
separating them
and consequently no three blocks have a point in common.
Furthermore, without loss of generality, we may assume that
the Gaussian curvature of every
$F_v$
is constant,
$K\equiv -1$,
and thus
$F_v$
can be viewed as a closed convex subset in
$\hyp^2$
bounded by countable many geodesic lines.

By Lemma~\ref{Lem:reduction2}, it suffices to construct a 
subexponential partition
$(X,g,T)$
of rank 1. Actually, there are several ways to do that,
for example, there is such an
$(X,g,T)$
with
$T=\R$
and 
$g^{-1}(t)$, $t\in T$,
homeomorphic to
$\R^2$
and quasi-isometric to a subset of a {\it template}
introduced in \cite{CK}. This easily implies that
each
$g^{-1}(t)$
has a subexponential growth rate, however, to show that it
has a bounded distortion in
$X$
is not that easy and it requires a large piece of machinery
from \cite{CK}. Apparently, this difficulty is rooted in the
fact that the target space
$T=\R$
is in a sense simplest possible.

We shall construct another subexponential partition
$(X,g,T)$
of rank 1 with
$T$
much more complicated for which it is not immediately
obvious even that
$\dim T=1$.
However, it is easy to establish the subexponential property of
$(X,g,T)$
as well as the fact that
$\rank(X,g,T)=1$.

To start with, we take a block
$X_v=F_v\times\R$,
fix a point
$p\in F_v$
and
$\rho>0$
with
$\dist(p,\d F_v)>2\rho$
and consider the metric projection
$\pr_v:F_v\sm B_{2\rho}(p)\to S_{2\rho}(p)$.
Then
$J_w=\pr_v(\d_wF_v)$
is an open interval in the circle
$S_{2\rho}(p)$
for every connected component
$\d_wF_v\sub\d F_v$.

Next, we note that
$I_w=\pr_v^{-1}(\ov J_w)\sub F_v$
has an subexponential growth rate because its area
is finite,
$\Area(I_w)<\pi$
by Gauss-Bonnet. Now, we collaps every
$\ov J_w$
to a point,
$\eta:S_{2\rho}(p)\to S$,
and identify the resulting factor-space 
$S=S_{2\rho}(p)/\sim$
with the circle
$S=S_{\rho}(p)\sub F_v$.
Furthermore, we may choose the identification in such a way
that the geodesic segment
$x\ov x$ 
lies in
$\ov B_{2\rho}(p)\sm B_{\rho}(p)$
for every
$x\in S_{2\rho}(p)$
and its image
$\ov x\in S$,
and different such segments have no common interior
points. Thus we can extend
$\eta$
to the closed annulus
$\ov B_{2\rho}(p)\sm B_{\rho}(p)$
by
$\eta(x')=\ov x$
for each
$x'\in x\ov x$.
Then 
$(F_v\sm B_{\rho}(p),\eta\circ\pr_v,S)$
is a subexponential partition because
$(\eta\circ\pr_v)^{-1}(x)$
is either a geodesic ray or some subset
$I_w$.
Note that there is only countable many points
$x\in S$
for which the last case occurs. Any such point
is said to be {\it distinguished}. Clearly, the set of
the distinguished points is a countable dense subset in
$S$.

Though the partition
$(F_v\sm B_{\rho}(p),\eta\circ\pr_v,S)$
has rank 1, we need to make further identifications
to extend it on the whole
$F_v$.
To this end, we fix a diameter
$T_v\sub\ov B_{\rho}(p)$
and take the metric projection
$\ov B_{\rho}(p)\to T_v$.
To avoid unnecessary complifications, we slightly
perturb it in such a way that any two distinguished 
points of
$S=\d\ov B_{\rho}(p)$
are mapped into distinct points and images of 
at most two points of
$S$
coincide. The images of the distinguised points 
are called distinguished too.
Let
$\zeta:\ov B_{\rho}(p)\to T_v$
be the perturbed map. Then 
$(F_v,\zeta\circ\eta\circ\pr_v,T_v)$
is a subexponential partition of rank 1. From now on,
we consider
$T_v$
as an abstract space and not as a subset of
$F_v$.
Let
$D_v\sub T_v$
be the subset of all distinguished points. Then
$\zeta\circ\eta\circ\pr_v$
maps bijectively the set
$\{\d_wF_v\}$
of the boundary components of
$F_v$
onto
$D_v$.

Coming back to the block
$X_v=F_v\times\R$,
we compose the projection on the first factor with
$\zeta\circ\eta\circ\pr_v$
and denote the resulting map by
$g_v:X_v\to T_v$.
Then, obviously,
$(X_v,g_v,T_v)$
is a subexponential partition of rank 1, and 
$g_v$
maps the set of the boundary components of
$X_v$
bijectively onto
$D_v$.

Finally, let
$\wt T=\sqcup_vT_v$
be the disjoint union. We define
$T:=\wt T/\sim$,
where distinct
$t\in T_v$, $t'\in T_{v'}$
are equivalent if and only if
$t\in D_v$, $t'\in D_{v'}$,
the blocks
$X_v$, $X_{v'}$
are different and adjacent along a common boundary component
which corresponds to
$t$
and
$t'$.
Then
$\wt g=\sqcup_vg_v:\sqcup_vX_v\to\sqcup_vT_v$
factors to a continuous
$g:X\to T$,
and the partition
$(X,g,T)$
is subexponential: for every distinguished
$t\in T$,
the set
$g^{-1}(t)$
sits in two adjacent blocks and it is
(outside of a 
$\{\text{compact set}\}\times\R$) 
the union of a geodesic half-plane and a subexponential
$I_w\times\R$
in each of the blocks. Otherwise,
$g^{-1}(t)$
is (again outside of a 
$\{\text{compact set}\}\times\R$) 
the union of two geodesic half-planes sitting 
in one and the same block. Each compact set
$K\sub X$
intersects only finitely many blocks
$X_v$,
hence its image
$g(K)\sub T$
intersects only finitely many segments
$T_v$
and consequently
$\dim g(K)\le 1$.
Thus
$\rank(X,g,T)=1$.
This completes the proof of Theorem~\ref{Thm:graphmfd}.


\bigskip
\noindent
Sergei Buyalo, Steklov Institute of Mathematics, Fontanka 27,
191011, St. Petersburg, Russia\\
{\tt buyalo@pdmi.ras.ru}

\medskip
\noindent
Viktor Schroeder, Institut f\"ur Mathematik, Universit\"at Z\"urich,
Winterthurer Strasse 190, CH-8057 Z\"urich, Switzerland\\
{\tt vschroed@math.unizh.ch}

\end{document}